\documentclass{article}
\PassOptionsToPackage{numbers,sort,compress}{natbib}
\usepackage{arxiv}    
\usepackage{natbib}
\usepackage{amsfonts}
\usepackage{amsmath,bm,mathtools,mathrsfs} 
\usepackage{amssymb}
\usepackage{amsthm}

\newtheorem{theorem}{Theorem}
\newtheorem{prop}[theorem]{Proposition}
\newtheorem{corr}[theorem]{Corollary}
\newtheorem{rmk}[theorem]{Remark}

\newcommand{\norm}[1]{\|{#1}\|}
\newcommand{\argmin}{\operatornamewithlimits{argmin}}

\begin{document} 
\title{ A Rate of Convergence for Two-Block Coordinate Descent}
\author{\name  Seyyed Mohammad Rouzban \\
\addr{University of Tehran} \\
\name Reshad Hosseini \email {reshad.hosseini@ut.ac.ir}\\
\addr{University of Tehran}\\
}
\maketitle

\vspace{2.5ex}
\begin{abstract}
This short report gives a non-asymptotic rate of convergence proof for solving a two-block coordinate descent problem. This non-asymptotic proof is a simple result that can be derived easily from available results in the literature. We give the results in this report because in this general form we have not seen being stated in the literature.
\end{abstract}

%%%%%%%%%%%%%%%%%%%%%%%%%
\section{Introduction}

In this short report, we show the global convergence guarantee for a class of two-block \textit{coordinate descent method} (BCD). We give a theorem for the global convergence rate of the proposed method. Our results can be derived easily from available results~\citep{nesterov2012efficiency,beck2013convergence} but we have not seen these results being stated clearly in those papers and other related literatures. Because of its wide applicability and importance, we give the results here.

%%%%%%%%%%%%%%%%%%%%%%%%%
\section{A Convergence Theorem for Coordinate Descent Method}
\label{sec:convergence}
In this part, we give the following convergence theorem for two-block coordinate descent. As a corollary, we show that the coordinate descent method, where we find the solution with respect to one block analytically.
\begin{theorem}\label{theorem.cond}
	Consider, we use a block coordinate descent algorithm for minimizing the differentiable function $f(x,y)$ with two blocks of variables $x$ and $y$. The block coordinate descent algorithm decreased the cost function with respect to block $y$ and finds a stationary point with respect to this block, i.e., $\nabla_y f(x_t,y_t)=0$. In the algorithm, the solution for the block $x$ satisfies the following condition
	\begin{equation}
	f(x_t,y_t)-f(x_{t+1},y_t)\ge \frac{1}{2E_{x_t}}\norm{\nabla_x f(x_t,y_t)}^2.\label{theorem.cond1}
	\end{equation}
	Then,
	\begin{equation}
	%\begin{split}
	\sum_{t=0}^{T-1}\frac{1}{E_{x_t}}\norm{\nabla f(x_t,y_t)}^2
	\le f(x_{0},y_{0})-f(x_{T},y_{T}).\label{lemma.conclusion}
	%\end{split}
	\end{equation}
	Furthermore, assume $f(x_t,y_t)$ is bounded from below and $E_{x_t}\geq E_{x}>0$ for all $t$ then every limit point is a stationary point. The rate of convergence to the stationary point is $O(1/\sqrt{T})$, where $T$ is the number of iterations.
\end{theorem}
\begin{proof}
	From~\eqref{theorem.cond1}, $\nabla_y f(x_t,y_t)=0$, and $f(x_{t+1},y_{t+1})\leq f(x_{t+1},y_t)$,  we have 	
	\begin{equation}
	\label{eq.thm1no1}
	\frac{1}{2E_{x_t}}\norm{\nabla f(x_t,y_t)}^2 \leq f(x_t,y_t)-f(x_{t+1},y_{t+1}). 
	\end{equation}	
	The summation of inequality~\eqref{eq.thm1no1} over $t$ from 0 to $T-1$, yields~\eqref{lemma.conclusion}. The right side of equation~\eqref{lemma.conclusion} is bounded and therefore if 
	$T\rightarrow \infty$ then $\nabla f(x^*,y^*)\rightarrow0$, where $(x^*,y^*)$ is a limit point. From~\eqref{lemma.conclusion}, it is easy to see that
	\begin{equation}
	\min_{0\leq t \leq T-1}\norm{\nabla f(x_t,y_t)}^2
	\le E_x\frac{f(x_{0},y_{0})-f(x_{T},y_{T})}{T},\nonumber
	\end{equation}
	where $E_x = \max_{1\leq t \leq T} E_{x_t}$. Therefore, we obtain the rate of convergence $O(1/\sqrt{T})$.
\end{proof}
The condition~\eqref{theorem.cond1} is very easy to be satisfied. For example, it is enough  to use  gradient descent with fixed step-size, and the gradient for the block is Lipschitz continuous (see~\citep{nesterov2004introductory}).

\begin{prop}
\label{prop.lip}
Assume the gradient function $g_t(x) = f(x,y_t)$ is Lipschitz continuous, i.e.,
\begin{equation*}
\|\nabla_x  f(x',y_t) - \nabla_x  f(x,y_t)\| \leq L(y_t) \| x' - x \|.
\end{equation*}
For the gradient descent step $x_{t+1} = x_t - \frac{1}{L(y_t)} \nabla_x  f(x_t,y_t)$, the following condition holds:
\begin{equation*}
f(x_{t},y_t) -  f(x_{t+1},y_t) \geq \frac{1}{2L(y_t)} \| \nabla_x  f(x_t,y_t) \|^2.
\end{equation*}
\end{prop}
If the optimization with respect to second block is solved analytically, $x_{t+1} = \argmin_x f(x,y_t)$. Then, if the function satisfies certain properties, the condition~\eqref{theorem.cond1} is satisfied. The following corollary is straightforward application of the previous proposition.
\begin{corr}
Assume the gradient function $g_t(x)=f(x,y_t)$ is Lipschitz continuous. If for the $x$ block we have $x_{t+1} = \argmin_x f(x,y_t)$, then
\begin{equation*}
f(x_{t},y_t) -  f(x_{t+1},y_t) \geq \frac{1}{2L(y_t)} \| \nabla_x  f(x_t,y_t) \|^2.
\end{equation*}
\end{corr}
\begin{proof}
Defining $x' = x_t - \frac{1}{L(y_t)} \nabla_x  f(x_t,y_t)$ and using the result of Proposition~\ref{prop.lip}, together with the fact  that $f(x_{t+1},y_t)\leq f(x',y_t)$, the proof is immediate.
\end{proof}
The previous result shows that two-block coordinate-descent method has a good convergence behavior far from optimum when the gradient is large. The convergence of gradient-descent is similar but its constant can be much larger depending on the Lipschitz constant of the whole function, while here it depends only on the Lipschitz constant for one block. Another important advantage also observed empirically for coordinate-descent method in compare to gradient decent is that here for any region of data it behaves like a gradient descent with best step-size while in practice the best step-size for gradient descent method is not available. With the following remark we finish this section.
\begin{rmk}
Without additional structure on the objective function, it seems hard to get better constants in the bound obtained. If the function $f(x,y_t)$ is  Lipschitz continuous gradient. Then if the function is twice differentiable, it is equal to say $\nabla_{xx} f(x,y_t) \leq L(y_t)$~\cite{nesterov2004introductory}. For the following quadratic function, the bound becomes equality:
	\begin{equation*}
	f(x,y_t)= c(y_t)+\nabla_x  f(x_t,y_t)^T (x-x_{t})+\tfrac{L(y_t)}{2} \| x-x_{t}\|^2.
	\end{equation*}
Therefore with the common assumption that assumes a Lipschitz constant for the set containing the whole iterations, it is hard to obtain better constant for the convergence rate.
\end{rmk}

\bibliographystyle{abbrvnat}
\bibliography{mybibfile}
\end{document}